\documentclass[12pt]{article}
\usepackage{epic,eepic}
\usepackage[dvips]{epsfig}
\usepackage{color}
\usepackage{amsmath, amscd, amssymb}
\usepackage{amssymb}
\def\draw #1 by #2 (#3){
  \vbox to #2{
    \hrule width #1 height 0pt depth 0pt
    \vfill
    \special{picture #3}
    }
  }
\def\scaleddraw #1 by #2 (#3 scaled #4){{
  \dimen0=#1 \dimen1=#2
  \divide\dimen0 by 1000 \multiply\dimen0 by #4
  \divide\dimen1 by 1000 \multiply\dimen1 by #4
  \draw \dimen0 by \dimen1 (#3 scaled #4)}
  }

\begin{document}
\renewcommand{\labelenumi}{\theenumi}
\newcommand{\qed}{\mbox{\raisebox{0.7ex}{\fbox{}}}}
\newtheorem{theorem}{Theorem}
\newtheorem{example}{Example}
\newtheorem{problem}[theorem]{Problem}
\newtheorem{defin}[theorem]{Definition}
\newtheorem{lemma}[theorem]{Lemma}
\newtheorem{corollary}[theorem]{Corollary}
\newtheorem{nt}{Note}
\newtheorem{proposition}[theorem]{Proposition}
\renewcommand{\thent}{}
\newenvironment{pf}{\medskip\noindent{\textbf{Proof}:  \hspace*{-.4cm}}\enspace}{\hfill \qed \medskip \newline}
\newenvironment{pft1}{\medskip\noindent{\textbf{Proof of Theorem~\ref{1}}:  \hspace*{-.4cm}}\enspace}{\hfill \qed \medskip \newline}
\newenvironment{pft8}{\medskip\noindent{\textbf{Proof of Theorem~\ref{8}}:  \hspace*{-.4cm}}\enspace}{\hfill \qed \medskip \newline}
\newenvironment{pft9}{\medskip\noindent{\textbf{Proof of Theorem~\ref{9}}:  \hspace*{-.4cm}}\enspace}{\hfill \qed \medskip \newline}
\newenvironment{defn}{\begin{defin}\em}{\end{defin}}{\vspace{-0.5cm}}
\newenvironment{lem}{\begin{lemma}\em}{\end{lemma}}{\vspace{-0.5cm}}
\newenvironment{cor}{\begin{corollary}\em}{\end{corollary}}{\vspace{-0.5cm}}
\newenvironment{thm}{\begin{theorem} \em}{\end{theorem}}{\vspace{-0.5cm}}
\newenvironment{pbm}{\begin{problem} \em}{\end{problem}}{\vspace{-0.5cm}}
\newenvironment{note}{\begin{nt} \em}{\end{nt}}{\vspace{-0.5cm}}
\newenvironment{exa}{\begin{example} \em}{\end{example}}{\vspace{-0.5cm}}
\newenvironment{pro}{\begin{proposition} \em}{\end{proposition}}{\vspace{-0.5cm}}

\setlength{\unitlength}{12pt}
\newcommand{\comb}[2]{\mbox{$\left(\!\!\begin{array}{c}
            {#1} \\[-0.5ex] {#2} \end{array}\!\!\right)$}}
\renewcommand{\labelenumi}{(\theenumi)}
\renewcommand{\b}{\beta}
\newcounter{myfig}
\newcounter{mytab}
\def\mod{\hbox{\rm mod }}
\def\scaleddraw #1 by #2 (#3 scaled #4){{
  \dimen0=#1 \dimen1=#2
  \divide\dimen0 by 1000 \multiply\dimen0 by #4
  \divide\dimen1 by 1000 \multiply\dimen1 by #4
  \draw \dimen0 by \dimen1 (#3 scaled #4)}
  }
\newcommand{\Aut}{\mbox{\rm Aut}}
\newcommand{\w}{\omega}
\def\r{\rho}
\newcommand{\DbF}{D \times^{\phi} F}
\newcommand{\autF}{{\tiny\Aut{\scriptscriptstyle(\!F\!)}}}
\def\Cay{\mbox{\rm Cay}}
\def\a{\alpha}
\newcommand{\C}[1]{\mathcal #1}
\newcommand{\B}[1]{\mathbb #1}
\newcommand{\F}[1]{\mathfrak #1}
\title{A relationship between the diameter and the intersection number $c_2$ for a distance-regular graph}

 \author{Jack H. Koolen \  and \  Jongyook Park\\
{\small {\tt koolen@postech.ac.kr} ~~
{\tt jongyook@postech.ac.kr}}\\
{\footnotesize{Department of Mathematics,  POSTECH, Pohang 790-785, South Korea}}\\
}

\date{\today}

\maketitle

\begin{abstract}
In this paper we will look at the relationship between the intersection number $c_2$ and its diameter for a distance-regular graph. And also, we give some tools to show that a distance-regular graph with large $c_2$ is bipartite, and a tool to show that if $k_D$ is too small then the distance-regular graph has to be antipodal.

\end{abstract}
\section{Introduction}

In this paper we will look at the relationship between the intersection number $c_2$ and its diameter for a distance-regular graph. (For definitions see next section.)

Let us first start with diameter three. A distance-regular graph $\Gamma$ with diameter three and valency $k$, can have $c_2 = k-1$. This occurs exactly when $\Gamma$ is the $K_{k+1, k+1}$ minus a perfect matching. But if $\Gamma$ is not bipartite, then it is fairly straightforward to show that $c_2 \leq \frac{2}{3}k$. (In Theorem~\ref{5} below we see that one of the following holds: $\Gamma$ is bipartite, a Taylor graph or $c_2 \leq \frac{1}{2}k$.)

We will show a similar behavior for $c_2$ when the diameter is larger, that is non-bipartite distance-regular graphs have a significant smaller $c_2$ than bipartite distance-regular graphs with the same diameter in general.

First, we concentrate on the situation when $\Gamma$ has diameter at least four containing a quadrangle.

Recall that Terwilliger (see \cite[Corollary 5.2.2]{bcn}) showed that if $\Gamma$ is a distance-regular graph with diameter $D$ and valency $k$, having an induced quadrangle, then $a_1 + 2 \leq \frac{2}{D}k$. Our first result shows that a similar result holds for the intersection number $c_2$ in this case:

\begin{thm}\label{1}
Let $\Gamma$ be a distance-regular graph with valency $k$ and diameter $D\geq4$. If $\Gamma$ contains a quadrangle, then $c_2\leq\frac{2}{D}k$. And equality if and only if  $D=4$ and $\Gamma$ is a Hadamard graph or $D\geq5$ and $\Gamma$ is a $D$-cube.
\end{thm}

\textbf{Remarks}: ($i$) Note that $K_{k+1, k+1}$ minus a perfect matching shows that for diameter three it is not true.

($ii$) For $c_2 =1$, the above result is not true, as the Foster and Biggs-Smith graphs show (with $D \geq 7$). But, if $k \geq 3$ and $D \geq 4$ we expect that $c_2 \leq \frac{2}{D}k$ holds, with a finite number of exceptions.\\

For diameters $4, 5, 6, 7$ we improve Theorem~\ref{1} and show that a similar behavior like diameter three occurs.

\begin{thm}\label{8}
Let $\Gamma$ be a distance-regular graph with valency $k\geq3$ and diameter $D$. If $D\in \{4,5\}$, then either $c_2\leq\frac{1}{3}k$ or $\Gamma$ is one of the following:
\begin{description}
  \item[$(i)$] $D=4$ and $\Gamma$ is a Hadamard graph,
  \item[$(ii)$] $D=5$ and $\Gamma$ is the $5$-cube.
\end{description}
\end{thm}

\begin{thm}\label{9}
Let $\Gamma$ be a distance-regular graph with valency $k\geq3$ and diameter $D$. If $D\geq6$, then either $c_2\leq\frac{1}{4}k$ or one of the following holds:
\begin{description}
  \item[$(i)$] $D=6$ and $\Gamma$ is either the generalized dodecagon of order $(1,2)$ or the $6$-cube,
  \item[$(ii)$] $D=7$ and $\Gamma$ is either the Biggs-Smith graph or the $7$-cube,
  \item[$(iii)$] $D=8$ and $\Gamma$ is the Foster graph.
\end{description}
\end{thm}

In order to show the above results, we first give preliminaries and definitions in next section. In Subsection 3.1, we give some tools to show that a distance-regular graph with large $c_2$ is bipartite, and in Subsection 3.2 we show that if $k_D$ is too small then the distance-regular graph has to be antipodal. In Section 4, we give a proof of Theorem~\ref{1}, and in Section 5 we show Theorem~\ref{8} and Theorem~\ref{9}.

\section{Definitions and preliminaries}
All the graphs considered in this paper are finite, undirected and
simple (for unexplained terminology and more details, see \cite{bcn}). Suppose
that $\Gamma$ is a connected graph with vertex set $V(\Gamma)$ and edge set $E(\Gamma)$, where $E(\Gamma)$ consists of unordered pairs of two adjacent vertices. The {\em distance} $d_{\Gamma}(x,y)$ between
any two vertices $x$ and $y$ in a graph $\Gamma$
is the length of a shortest path connecting $x$ and $y$. If the graph $\Gamma$ is clear from the context, then we simply use $d(x,y)$. We define the {\em diameter} $D$
of $\Gamma$ as the maximum distance in $\Gamma$.  For a vertex $x \in V(\Gamma)$, define $\Gamma_i(x)$ to be the set of
vertices which are at distance precisely $i$ from $x~(0\le i\le D)$. In addition, define $\Gamma_{-1}(x) = \Gamma_{D+1}(x)
:= \emptyset$. We write $\Gamma(x)$ instead of $ \Gamma_1(x)$.

A connected graph $\Gamma$ with diameter $D$ is called {\em{distance-regular}} if there are integers $b_i, c_i$ $(0\leq i\leq D)$ such that for any two vertices $x, y \in V(\Gamma)$ with $d(x, y)=i$, there are precisely $c_i$ neighbors of $y$ in $\Gamma_{i-1}(x)$ and $b_i$ neighbors of $y$ in $\Gamma_{i+1}(x)$, where we define $b_D=c_0=0$. In particular, any distance-regular graph  is regular with valency $k := b_0$. Note that a  (non-complete) connected {\em strongly regular graph} is just a distance-regular graph with diameter two. We define $a_i := k-b_i-c_i$ for notational convenience.  Note that $a_i=\mid\Gamma(y)\cap\Gamma_i(x)\mid$ holds for any two vertices $x, y$ with $d(x, y)=i$ $(0\leq i\leq D).$

 For a distance-regular graph $\Gamma$ and a vertex $x\in V(\Gamma)$, we denote $k_i:=|\Gamma_i(x)|$ and $p^i_{jh}:=|\{w| w\in\Gamma_j(x)\cap\Gamma_h(y)\}|$ for any  $y\in\Gamma_i(x)$. It is easy to see that $k_i = \frac{b_0 b_1 \cdots b_{i-1}}{c_1 c_2 \cdots c_i}$ and hence it does not depend on $x$.  The numbers $a_i$, $b_{i-1}$ and $c_i$ $(1\leq i\leq D)$ are called the {\em{intersection~numbers}}, and the array $\{b_0,b_1,\cdots,b_{D-1};c_1,c_2,\cdots,c_D\}$ is called the {\em{intersection~array}} of $\Gamma$. A distance-regular graph with intersection array $\{k,\mu,1;1,\mu,k\}$ is called a {\em Taylor graph}.

 Suppose that $\Gamma$ is a distance-regular graph with valency $k \geq 2$ and diameter $D \geq 2$, and let $A_i$ be the matrix of $\Gamma$ such that the rows and the columns of $A_i$ are indexed by the vertices of $\Gamma$ and the ($x, y$)-entry is $1$ whenever $x$ and $y$ are at distance $i$ and $0$ otherwise. We will denote the adjacency matrix of $\Gamma$ as $A$ instead of $A_1$. The eigenvalues of the graph $\Gamma$ are the eigenvalues of $A$.

Some standard properties of the intersection numbers are collected in the following lemma.

\begin{lem}(\cite[Proposition 4.1.6]{bcn})\label{pre}{\ \\}
Let $\Gamma$ be a distance-regular graph with valency $k$ and diameter $D$. Then the
following holds:\\
$(i)$ $k=b_0> b_1\geq \cdots \geq b_{D-1}~;$\\
$(ii)$ $1=c_1\leq c_2\leq \cdots \leq c_{D}~;$\\
$(iii)$ $b_i\ge c_j$ \mbox{ if }$i+j\le D~.$
\end{lem}

Suppose that $\Gamma$ is a distance-regular graph with valency $k\ge 2$ and diameter $D\ge 1$. Then $\Gamma$ has exactly $D+1$ distinct eigenvalues, namely $k=\theta_0>\theta_1>\cdots>\theta_D$ (\cite[p.128]{bcn}).

Recall that a {\em clique} of a graph is a set of mutually adjacent vertices and that a {\em co-clique} of a graph is a set of vertices with no edges.  A clique $\mathcal{C}$ of a distance-regular graph with valency $k$, diameter $D\geq2$ and smallest eigenvalue $\theta_D$, is called a {\em Delsarte clique} if $\mathcal{C}$ contains exactly $1-\frac{k}{\theta_D}$ vertices.

A graph $\Gamma$ is called {\em bipartite} if it has no odd cycle. (If $\Gamma$ is a distance-regular graph with diameter $D$ and bipartite, then $a_1=a_2=\ldots=a_D=0$.)  An {\em antipodal} graph is a connected graph $\Gamma$ with diameter $D>1$ for which being at distance 0 or $D$ is an equivalence relation. If, moreover, all equivalence classes have the same size $r$, then $\Gamma$ is also called an {\em antipodal $r$-cover}. \\

Recall the following results.

\begin{thm}(cf. \cite[Proposition 5]{kp2})\label{5}
Let $\Gamma$ be a distance-regular graph with valency $k$ and diameter $D$. If $D=3$, then one of the following holds:
\begin{description}
  \item[$(i)$] $c_2\leq\frac{1}{2}k$, $b_2\leq\frac{1}{2}k_3$ and $c_3\leq\frac{1}{2}k_2$,
  \item[$(ii)$] $\Gamma$ is bipartite,
  \item[$(iii)$] $\Gamma$ is a Taylor graph.
\end{description}
\end{thm}

\begin{lem}\label{del}(cf.\cite[Proposition 4.4.6]{bcn}) Let $\Gamma$ be a distance-regular graph with diameter $D\geq2$ and distinct eigenvalues $k=\theta_0>\theta_1>\cdots>\theta_D$. Then:
\begin{description}
  \item[$(i)$] The number of vertices of a clique $\mathcal{C}$ in $\Gamma$ is bounded by $$\mid\mathcal{C}\mid\leq 1-\frac{k}{\theta_D}.$$ Moreover, equality holds if and only if the clique $\mathcal{C}$ is a Delsarte clique.
  \item[$(ii)$] If $\Gamma$ contains a nonempty induced complete bipartite subgraph $K_{s,t}$, then $$\frac{2st}{s+t}\leq \frac{b_1}{\theta_1+1}+1.$$
\end{description}
\end{lem}

\begin{thm}\label{inter}(cf.\cite{heam})   Let $m \leq n$ be two positive integers. Let
$A$ be an $n\times n$ matrix, that is similar to a (real) symmetric matrix, and let
$B$ be a principal $m \times m$ submatrix of $A$. Then, for $i=1,\ldots , m$, $$\theta_{n-m+i}(A)\leq \theta_i(B)\leq \theta_i(A)$$
holds, where $A$ has eigenvalues $\theta_1(A) \geq \theta_2(A) \geq \ldots\geq \theta_n(A)$ and B has eigenvalues
$\theta_1(B) \geq \theta_2(B) \geq \ldots \geq \theta_m(B)$.\\
\end{thm}

\begin{lem}\label{4}
Let $\Gamma$ be a distance-regular graph with valency $k$ and diameter $D$. If $a_1\neq0$, then $c_{\lfloor\frac{D}{2}\rfloor}\leq\frac{1}{3}k$
\end{lem}
\begin{pf}
This follows immediately from \cite[Proposition 5.5.6]{bcn}.
\end{pf}

\section{Some preliminary results}

In this section, we give some preliminary results which are helpful to prove our results in section 4 and section 5. We first give some tools to show that a distance-regular graph with large $c_2$ is bipartite, and then we show that if $k_D$ is too small then the distance-regular graph has to be antipodal.

\subsection{Tools to show bipartiteness}

The following lemma is useful to check whether a distance-regular graph is bipartite, from a condition of the intersection numbers $c_i$.

\begin{lem}\label{3}
Let $\Gamma$ be a distance-regular graph with valency $k\geq3$ and diameter $D\geq4$. For an integer $i\in\{2,3,\ldots,D\}$, if $a_i\neq0$ and $2c_i+c_{D-i}>k$, then $a_{i-1}\neq0$.
\end{lem}

\begin{pf}
Let $x,y$ and $z$ be vertices such that $d(x,y)=i$, $d(x,z)=D-i$ and $d(y,z)=D$. As $a_i\neq0$, there exists $y'$ which is adjacent to $y$ and at distance $i$ from $x$. Then $\Gamma(x)\cap\Gamma_{D-i-1}(z)$ has cardinality $c_{D-i}$ and is disjoint from $\Gamma(x)\cap(\Gamma_{i-1}(y)\cup\Gamma_{i-1}(y'))$. Hence, $k=|\Gamma(x)|\geq|\Gamma(x)\cap\Gamma_{D-i-1}(z)|+(|\Gamma(x)\cap\Gamma_{i-1}(y)|+|\Gamma(x)\cap\Gamma_{i-1}(y')|-|\Gamma(x)\cap(\Gamma_{i-1}(y)\cap\Gamma_{i-1}(y'))|)\geq c_{D-i}+2c_i-p^1_{i-1i-1}=c_{D-i}+2c_i-\frac{k_{i-1}}{k}a_{i-1}$. This shows the lemma.
\end{pf}

As a consequence of previous two lemmas, we have the following lemma.

\begin{lem}\label{6}
Let $\Gamma$ be a distance-regular graph with valency $k\geq3$ and diameter $D\geq4$. If $c_2>\frac{1}{3}k$, then $D\leq5$ and $\Gamma$ is bipartite.
\end{lem}

\begin{pf}
As $c_2>\frac{1}{3}k$ and $D\geq4$, we know $a_1=0$ by Lemma~\ref{4}. Hence, $a_2=0$ by Lemma~\ref{3}. Now, we have $c_3\geq\frac{3}{2}c_2>\frac{1}{2}k$ by \cite[Theorem5.4.1]{bcn}. Then, we have $a_i=0$ for $i\geq3$ by Lemma~\ref{3}. So, $\Gamma$ is bipartite. Now, $c_3>\frac{1}{2}k$ implies $c_3>b_3$ and hence $D\leq5$.
\end{pf}

The next lemma shows that if $\Gamma$ is a bipartite distance-regular graph with even diameter, then the intersection number $c_2$ divides the valency $k$.

\begin{lem}\label{2}
Let $\Gamma$ be a distance-regular graph with valency $k\geq3$ and diameter $D=2t$ for some integer $t\geq2$. If $\Gamma$ is bipartite, then $k_2=\alpha(k-1)$ for some integer $\alpha$.
\end{lem}
\begin{pf}
Let $\Gamma^{\frac{1}{2}}$ be the halved graph of $\Gamma$, and let $x$ be a vertex of $\Gamma$. We may assume that $\Gamma(x)\cup\Gamma_3(x)\cup\cdots\cup\Gamma_{D-1}(x)$ is the vertex set of $\Gamma^{\frac{1}{2}}$. Here note that the halved graph $\Gamma^{\frac{1}{2}}$ has valency $k_2$ and diameter $t$. As $\Gamma$ is bipartite, for any two vertices  $y,z\in\Gamma(x)$, $d_{\Gamma}(y,z)=2$. i.e. $\Gamma^{\frac{1}{2}}$ contains a clique $C$ with $k=\mid\Gamma(x)\mid$ vertices. As $2d_{\Gamma^{\frac{1}{2}}}(C,w)=d_{\Gamma}(C,w)\leq D-2=2t-2$ for any vertex $w$ of $\Gamma^{\frac{1}{2}}$, one can easily see that $C$ is a completely regular code with covering radius $t-1$ (see, \cite[Chapter 11.1]{bcn}), and hence $C$ is a Delsarte clique of $\Gamma^{\frac{1}{2}}$. So, $k=\mid C\mid=1-\frac{k_2}{\theta_{min}(\Gamma^{\frac{1}{2}})}$, where ${\theta_{min}(\Gamma^{\frac{1}{2}})}$ is the smallest eigenvalue of $\Gamma^{\frac{1}{2}}$. As $\theta_{min}(\Gamma^{\frac{1}{2}})$ is an algebraic integer, $\theta_{min}(\Gamma^{\frac{1}{2}})$ should be an integer.
\end{pf}

\subsection{A tool to show antipodality}
The following theorem is helpful to check whether a distance-regular graph with small $k_D$ is antipodal.

\begin{thm}\label{7}
Let $\Gamma$ be a distance-regular graph with valency $k\geq3$ and diameter $D\geq3$. If $a_D=0$ and $k_{D-1}<2k$, then one of the following holds:
\begin{description}
  \item[$(i)$]  $k_D=1$ and $\Gamma$ is an antipodal $2$-cover;
  \item[$(ii)$] $D=3$ and $\Gamma$ is bipartite.
\end{description}
\end{thm}

\begin{pf}
If $k_D=1$, then the graph $\Gamma$ is an antipodal $2$-cover. So, we may assume $k_D\geq2$. If $D=3$, then by Theorem~\ref{5}, $k_{D-1}\geq2c_D=2k$ or $\Gamma$ is bipartite. So, the theorem is true for $D=3$. Hence, we may assume $D\geq4$. As $b_2\geq c_2$(Lemma~\ref{pre}), we have $c_2\leq\frac{1}{2}k$. As $a_D=0$, we have $p^D_{D2}=\frac{a_D(a_D-1-a_1)+c_D(b_{D-1}-1)}{c_2}=\frac{k(b_{D-1}-1)}{c_2}$, and this implies $b_{D-1}=\frac{c_2p^D_{D2}}{k}+1\leq\frac{c_2(k_D-1)}{k}+1\leq\frac{1}{2}(k_D+1)$, as $c_2\leq\frac{1}{2}k$. But $k_{D-1}<2k$ implies $b_{D-1}>\frac{1}{2}k_D$. So, $b_{D-1}$ should be equal to $\frac{1}{2}(k_D+1)$ as $k_D$ is an integer, and this implies equality in the previous inequality. i.e. $c_2=\frac{1}{2}k$. Then by Lemma~\ref{6}, $\Gamma$ is bipartite and $D\leq4$, as $c_3\geq\frac{3}{2}c_2=\frac{3}{4}k$ and $b_2=\frac{1}{2}k$.\\

For $D=4$, we put $b_3=\alpha$, whence $k_4=2\alpha-1$ and $c_3=k-\alpha$. As $k_4$ is odd and more than one, $k_4\geq3$ and hence $\alpha\geq2$. As $kb_1=c_2k_2$ and $k_2b_2=c_3k_3$, we have $k_2=2(k-1)$ and $k_3=\frac{k-1}{k-\alpha}k$. Then, we find $2\alpha-1=k_4=\frac{b_3}{c_4}k_3=\frac{k-1}{k-\alpha}\alpha$, and this gives $(2-\frac{k-1}{k-\alpha})\alpha=1$. As $k-\alpha=c_3\geq\frac{3}{4}k$, we know $\frac{k-1}{k-\alpha}\leq\frac{k-1}{3k/4}<\frac{4}{3}$, and this implies $1=(2-\frac{k-1}{k-\alpha})\alpha>\frac{2}{3}\alpha\geq\frac{4}{3}$. This is a contradiction. This shows the theorem.\\

\end{pf}

\section{Proof of Theorem~\ref{1}}
In this section we give a proof of Theorem~\ref{1}. Before showing Theorem~\ref{1}, we first show Proposition~\ref{10}. Then we use this proposition to prove Theorem~\ref{1}.

\begin{pro}\label{10}
Let $\Gamma$ be a distance-regular graph with valency $k$ and diameter $D$. Let $t\geq2$ be an integer. If $D\geq2t$, $c_2>\frac{1}{t+1}k$ and $c_i>c_{i-1}$ holds for any $i\in\{2,3,\ldots,D\}$, then $D\in \{2t,2t+1\}$ and one of the following holds:
\begin{description}
  \item[$(i)$] $D=4$ and $\Gamma$ is a Hadamard graph,
  \item[$(ii)$] $D\geq5$ and $\Gamma$ is a $D$-cube.
\end{description}
\end{pro}

\begin{pf}
As the sequence $(c_i)_{i=1,\ldots,D}$ is a strictly increasing sequence, we know that $c_i\geq\frac{i}{2}c_2$ holds for $2\leq i\leq D$ by \cite[Theorem 5.4.1]{bcn} and \cite[Proposition 1 (ii)]{koolen}, and this implies $c_{\lfloor\frac{D}{2}\rfloor}\geq c_t\geq\frac{t}{2}c_2>\frac{t}{2(t+1)}k\geq\frac{1}{3}k$, as $t\geq2$. By Lemma~\ref{4}, we have $a_1=0$. Also, $2c_i+c_{D-i}\geq ic_2+\frac{D-i}{2}c_2\geq ic_2+\frac{2t-i}{2}c_2>\frac{2t+i}{2}\times\frac{1}{t+1}k\geq k$ holds for any $i\in\{2,3,\ldots,D\}$. i.e. $a_2=a_3=\cdots= a_D=0$ by Lemma~\ref{3}. So, the graph $\Gamma$ is bipartite. As $c_{t+1}>\frac{1}{2}k$, $\Gamma$ has diameter at most $2t+1$ by Lemma~\ref{pre}. i.e. $D\in\{2t,2t+1\}$.\\

 Here note that if $D=2t$, then $k_2=\frac{k(k-1)}{c_2}<(t+1)(k-1)$ should divide $k-1$ by Lemma~\ref{2}, and this shows $c_2\geq\frac{1}{t}k$. But then $D=2t$, and $k-\frac{D-i}{2}c_2\geq k-c_{D-i}=b_{D-i}\geq c_i\geq\frac{i}{2}c_2$ ($2\leq i\leq D$) implies that $c_2\leq\frac{1}{t}k$. This means that $c_2=\frac{1}{t}k$ and $c_i=b_{D-i}=\frac{i}{2}c_2=\frac{i}{2t}k$ holds for $i\in\{2,3,\ldots D-2\}$ (where $c_1=1$ and $b_1=k-1$). Then one can easily see that $k_{D-1}<2k$ and this shows that the graph $\Gamma$ is an antipodal $2$-cover by Theorem~\ref{7}.\\

If $t=2$, then $D\in\{4,5\}$. For $D=4$, one can easily show that $\Gamma$ has an intersection array $\{k,k-1,\frac{1}{2}k,1;1,\frac{1}{2}k,k-1,k\}$, as $c_2=\frac{1}{2}k$ and $\Gamma$ is antipodal and bipartite. So, $\Gamma$ is a Hadamard graph.\\
Now we assume $D=5$, and consider the halved graph $\Gamma^{\frac{1}{2}}$ of $\Gamma$ to show that the graph $\Gamma$ has at most $132$ vertices. And then we can check the feasible intersection arrays in \cite[p.418]{bcn} to show that the graph $\Gamma$ is the $5$-cube. Clearly, the halved graph $\Gamma^{\frac{1}{2}}$ of $\Gamma$ is a strongly regular graph with valency $k_2$ and diameter $2$. As $c_2>\frac{1}{3}k$, $c_3>\frac{1}{2}k$ and $c_4>\frac{2}{3}k$, $\Gamma$ has at most $12k$ vertices, as $k_2<3(k-1)$, $k_3<4(k-1)$, $k_4<3(k-1)$ and $k_5<k-1$.\\
For a fixed vertex $x$ of $\Gamma$, we may assume that $\Gamma^{\frac{1}{2}}$ has vertex set $\Gamma(x)\cup\Gamma_3(x)\cup\Gamma_5(x)$. Then, the set $\Gamma(x)$ is a clique of $\Gamma^{\frac{1}{2}}$, and this implies that the smallest eigenvalue $\theta_{\rm min}(\Gamma^{\frac{1}{2}})$ of $\Gamma^{\frac{1}{2}}$ is bigger than $-3$ by Lemma~\ref{del} $(i)$. So, $\theta_{\rm min}(\Gamma^{\frac{1}{2}})=-2$ or $\Gamma^{\frac{1}{2}}$ is a conference graph, as  $\theta_{\rm min}(\Gamma^{\frac{1}{2}})=-1$ implies that $\Gamma^{\frac{1}{2}}$ is a complete graph.\\
If $\Gamma^{\frac{1}{2}}$ is a conference graph, then one can see that $\Gamma^{\frac{1}{2}}$ has at most $21$ vertices, as $\theta_{\rm min}(\Gamma^{\frac{1}{2}})>-3$. i.e. $\Gamma$ has at most $42$ vertices in this case.\\
So, we may assume $\theta_{\rm min}(\Gamma^{\frac{1}{2}})=-2$. As $\Gamma$ has at most $12k$ vertices, the halved graph $\Gamma^{\frac{1}{2}}$ has at most $6k$ vertices, and this implies that the number of vertices of $\Gamma^{\frac{1}{2}}$ is at most $6\times$(the maximum number of vertices of any clique in $\Gamma^{\frac{1}{2}}$), as $k=\mid\Gamma(x)\mid\leq$ the maximum number of vertices of any clique in $\Gamma^{\frac{1}{2}}$. As $\theta_{\rm min}(\Gamma^{\frac{1}{2}})=-2$, we know that $\Gamma^{\frac{1}{2}}$ is a triangular graph $T(n)$, a square grid $n\times n$, a complete multipartite graph $K_{n\times2}$, or one of the graphs of Petersen, Clebsch, Schlafli, Shrikhande, or Chang by \cite[Theorem 3.12.4]{bcn}.\\
If $\Gamma^{\frac{1}{2}}$ is a triangular graph $T(n)$, then $\Gamma^{\frac{1}{2}}$ has at most $66$ vertices, as $T(n)$ has $\frac{n(n-1)}{2}$ vertices and a maximum clique with $n-1$ vertices.\\
If $\Gamma^{\frac{1}{2}}$ is a square grid $n\times n$, then $\Gamma^{\frac{1}{2}}$ has at most $36$ vertices, as $n\times n$ has $n^2$ vertices and a maximum clique with $n$ vertices\\
If $\Gamma^{\frac{1}{2}}$ is a complete multipartite graph $K_{n\times2}$, then $\Gamma^{\frac{1}{2}}$ has at most $12$ vertices, as $K_{n\times 2}$ has $2n$ vertices and a maximum clique with $2$ vertices\\
If $\Gamma^{\frac{1}{2}}$ is one of the graphs of Petersen, Clebsch, Schlafli, Shrikhande, or Chang, then  $\Gamma^{\frac{1}{2}}$ has at most $28$ vertices (see, \cite[p.103-105]{bcn}).\\
So, the graph $\Gamma$ has at most $132$ vertices. Then we check the feasible intersection arrays in \cite[p.418]{bcn}, and we find that the only possible case is that the graph $\Gamma$ is the $5$-cube. Here note that the Clebsch graph is isomorphic to the halved $5$-cube.\\

 Now, we assume $t\geq3$. i.e. $D\geq6$. If $D=6$,  then $\Gamma$ is the $6$-cube by \cite[Corollary 5.8.3]{bcn}, as $\Gamma$ is bipartite and an antipodal $2$-cover. So, we may assume $D\geq7$. We first show that the second largest eigenvalue $\theta_1$ of $\Gamma$ is at least $\frac{1}{2}k$, and then we show the graph $\Gamma$ is the $D$-cube by showing $c_2=2$.\\
 Let $x$ and $y$ be vertices of $\Gamma$ such that they are at distance $D$ in $\Gamma$. Let $\Gamma_{0}^{t-1}$  be the induced subgraph of $\Gamma$ such that it has vertex set ${\displaystyle\bigcup^{t-1}_{i=0}\Gamma_i(x)}\cup\displaystyle\bigcup^{t-1}_{j=0}\Gamma_j(y)$. Here note that $\Gamma_{0}^{t-1}$ of $\Gamma$ consists of two disjoint isomorphic components, and hence the second largest eigenvalue $\theta_1$ of $\Gamma$ is at least the largest eigenvalue of a component of $\Gamma_{0}^{t-1}$ by Theorem~\ref{inter}.\\
 Let us consider the right lower $3\times3$ principal submatrix $P_t$ of $Q(\Gamma_0^{t-1})$. As $P_t=\left[
                 \begin{array}{ccc}
                   0 & b_{t-3} & 0 \\
                   c_{t-2} & 0 & b_{t-2} \\
                   0 & c_{t-1} & 0 \\
                 \end{array}
               \right]$ has eigenvalues $0, \pm\sqrt{c_{t-1}b_{t-2}+c_{t-2}b_{t-3}}$, the second largest eigenvalue $\theta_1$ of $\Gamma$ is at least $\sqrt{c_{t-1}b_{t-2}+c_{t-2}b_{t-3}}$.\\
 If $t=3$, then $\theta_1\geq\sqrt{c_{t-1}b_{t-2}+c_{t-2}b_{t-3}}=\sqrt{c_2(k-1)+k}>\sqrt{\frac{1}{4}k(k+3)}>\frac{1}{2}k$ holds, as $c_2>\frac{1}{t+1}k=\frac{1}{4}k$.\\
 If $t\geq4$, then  $\theta_1\geq\sqrt{c_{t-1}b_{t-2}+c_{t-2}b_{t-3}}\geq \sqrt{c_{t-1}c_{t+2}+c_{t-2}c_{t+3}}>\frac{1}{2}k$ holds, as $D\geq2t$ and $c_i\geq\frac{i}{2}c_2>\frac{i}{2(t+1)}k$ holds for $D\geq i\geq2$. So, the second largest eigenvalue of $\theta_1$ of $\Gamma$ is at least $\frac{1}{2}k$. \\
 Now, we show $c_2=2$. Take a vertex $z$ in $\Gamma_2(x)$, then the induced subgraph of $\Gamma$ on $\{x,z\}\cup(\Gamma(x)\cap\Gamma(z))$ is a complete bipartite $K_{2,c_2}$. By Lemma~\ref{del} $(ii)$, we have $\frac{2\times2\times c_2}{c_2+2}\leq\frac{b_1}{\theta_1+1}+1<3$. i.e. $c_2\in\{2,3,4,5\}$.\\
To show $c_2=2$, we first show that $2c_2>c_3$ holds. Clearly, the inequality is true for even diameter. So, we assume that the graph $\Gamma$ has odd diameter and $c_3\geq2c_2$ holds. Then one can show that $c_{2j+1}\geq (j+1)c_2$ holds for $j\in\{1,2,\ldots,t\}$ by \cite[Proposition 1 (ii)]{koolen}. i.e. $c_{2t+1}\geq(t+1)c_2>k$. This is a contradiction. So, $2c_2>c_3$ holds. Then, by \cite[Theorem 2.2]{cau}, we have $c_3-1-c_2(c_2-1)\geq-\frac{c_2(c_2-1)(c_2-2)^2}{2b_2}$, and this implies  $2c_2-1-c_2(c_2-1)>-\frac{c_2(c_2-1)(c_2-2)^2}{2b_2}$, as $2c_2>c_3$.\\
If $c_2\in\{3,4,5\}$, then the previous inequality shows $b_2\leq 2c_2$, which implies $D\leq6$ as $2c_2\leq c_4$. This contradicts $D\geq7$.\\
So, we find $c_2=2$ and this implies $k\leq 2t+1$, as $2=c_2>\frac{1}{t+1}k$.\\
If $D=2t$, then clearly the graph $\Gamma$ is the $D$-cube, as $c_i=b_{D-i}=\frac{i}{2}c_2=i$ for any $i\in\{1,\ldots,D\}$.\\
If $D=2t+1$, then one can easily see that $c_i=i$ holds for any $i\in\{1,\ldots,D\}$, as $c_2=2$, $c_{2t+1}=c_D=k\leq2t+1$ and $c_i>c_{i-1}$ holds for any $i\in\{2,\ldots,D\}$. So, $c_i=i$ for any $i\in\{1,\ldots,D\}$, and hence the graph $\Gamma$ is the $D$-cube, as $k=2t+1$ and $b_i=c_{D-i}$($0\leq i\leq D-1$).\\

This shows the proposition.
\end{pf}

Now, we give a proof of Theorem~\ref{1}

\begin{pft1}
As $\Gamma$ contains a quadrangle, we have $D\leq\frac{2}{a_1+2}k$ by \cite[Corollary 5.2.2]{bcn}. i.e. $a_1+2\leq\frac{2}{D}k$.\\
Hence if $c_2\leq a_1+2$, then clearly $c_2\leq\frac{2}{D}k$ holds. Also, $c_2=\frac{2}{D}k$ implies $c_2\geq a_1+2$. So, we may assume $c_2\geq a_1+3$. Then the sequence $(c_i)_{i=1,\ldots,D}$ is strictly increasing by \cite[Theorem 5.2.5]{bcn}, and hence  $c_i\geq\frac{i}{2}c_2$ holds for $2\leq i\leq D$ by \cite[Theorem 5.4.1]{bcn} and \cite[Proposition 1 (ii)]{koolen}. So, $c_2\leq \frac{2}{D}c_D\leq\frac{2}{D}k$ holds.\\

Now, we assume that equality holds. Then, $c_2\geq a_1+2>a_1$, and hence the sequence $(c_i)_{i=1,\ldots,D}$ is strictly increasing. This in turn shows that $\frac{i}{D}k=b_{D-i}=c_i$ hold for $i\in\{2,3,\ldots,D-2\}$. As a consequence we find $a_i=0$ for  $i\in\{2,3,\ldots,D-2\}$. As $a_2=0$, we have $a_1=0$ by \cite[Proposition 5.5.1]{bcn}, and \cite[Proposition 5.5.4]{bcn} implies $a_{D-1}=a_D=0$, as $a_{D-2}=0$ and $c_{D-2}>\frac{1}{2}k$. So, $\Gamma$ is bipartite. As $\frac{i}{D}k=b_{D-i}=c_i$ ($i\in\{2,3,\ldots,D-2\}$), $c_{D-1}\geq c_{D-2}=\frac{D-2}{D}k$ and $D\geq4$, we have $k_{D-1}=\frac{kb_1b_2\cdots b_{D-2}}{c_2c_3\cdots c_{D-1}}\leq\frac{D}{D-2}(k-1)<2k$. Then, by Theorem~\ref{7}, $\Gamma$ is an antipodal $2$-cover.\\

Let $t$ be the maximal integer such that $D\geq2t$. Then $D\in\{2t,2t+1\}$. For both cases, $\frac{2}{D}k=c_2>\frac{1}{t+1}k$ holds. So, the graph $\Gamma$ is either a Hadamard graph($D=4$) or the $D$-cube by Proposition~\ref{10}.\\

One can easily see that the converse is true.
\end{pft1}

\section{Proofs of Theorem~\ref{8} and Theorem~\ref{9}}
In this section, we give the proofs of Theorem~\ref{8} and Theorem~\ref{9}.

\begin{pft8}
Assume $D\geq4$ and $c_2>\frac{1}{3}k$, then by Lemma~\ref{6}, $D\leq5$ and $\Gamma$ is bipartite. Now $c_2>a_1$ implies that the sequence $(c_i)_{i=1,\ldots,D}$ is strictly increasing, and, by Proposition~\ref{10}, we see that the graph $\Gamma$ is a Hadamard graph or the $5$-cube.\\

\end{pft8}

\begin{pft9}
Suppose $c_2>\frac{1}{4}k$. Here note that if $k=3$, then the graph $\Gamma$ is the generalized dodecagon of order $(1,2)$, the Biggs-Smith graph or the Foster graph by \cite[Theorem 7.5.1]{bcn}(cf. \cite[Table 1]{bbs}). So, we may assume $k\geq4$, and hence $c_2\geq2$. As $c_2\geq2$, we have $c_3\geq\frac{3}{2}c_2\geq\frac{3}{8}k$ by \cite[Theorem 5.4.1]{bcn}. Then by Lemma~\ref{4}, we find $a_1=0$. As $c_2\geq2$ and $c_2>a_1$, by \cite[Theorem 5.2.5]{bcn} the sequence $(c_i)_{i=1,\ldots,D}$ is strictly increasing. By taking $t=3$ in Proposition~\ref{10}, the graph $\Gamma$ is either the $6$-cube or the $7$-cube.
\end{pft9}

\section{Acknowledgments}
J.H. Koolen was partially supported by the Basic Science Research Program
through the National Research Foundation of Korea(NRF) funded by the Ministry
of Education, Science and Technology (Grant number 2009-0089826).

\clearpage


\bigskip



\bigskip

\begin{thebibliography}{99}


\bibitem{bbs}
N.L. Biggs, A.G. Boshier and J. Shawe-Taylor, Cubic distance-regular graphs, J. London Math. Soc. 33 (1986) 385-394



\bibitem{bcn}
A.E. Brouwer, A.M. Cohen and A. Neumaier, Distance-Regular Graphs, Springer-Verlag, Berlin, 1989.



\bibitem{cau}
J.S. Caughman, IV, Intersection numbers of bipartite distance-regular graphs, Discrete Math. 163 (1997) 235-241



\bibitem{heam}
W.H. Haemers, Interlacing eigenvalues and graphs, Linear Algebra Appl. 226/228 (1995), 593-616.

\bibitem{koolen}
J.H. Koolen, On subgraphs in distance-regular graphs, J. Algebraic Combin. 1 (4) (1992), 353-362.


\bibitem{kp2}
J.H. Koolen, Jongyook Park, Distance-regular graphs with large $a_1$ or $c_2$, arXiv:1008.1209v1.






\end{thebibliography}
\end{document}